\documentclass[preprint,3p,11pt]{elsarticle}

 \usepackage{hyperref}
 \hypersetup{backref, colorlinks=true, linkcolor=blue}

\usepackage{stfloats}
\usepackage{subcaption}
\usepackage{placeins}

\usepackage{graphicx} 
\usepackage{booktabs} 
\usepackage{amssymb,amsmath}
\usepackage[ruled,vlined]{algorithm2e}
\usepackage{xcolor}

\begin{document}
\begin{center}
{\bf \Large Hybrid Two-Level  Transport Method with Solution Decomposition\\
	\vspace{0.1cm}
	 in Macro and Micro Components}
\end{center}
\begin{center}
Caleb A. Shaw and Dmitriy Y. Anistratov 
\end{center}
\begin{center}
{ \it Department of Nuclear Engineering,
North Carolina State University,
Raleigh, NC 27695 \\
 cashaw4@ncsu.edu, anistratov@ncsu.edu }
\end{center}

\begin{frontmatter}
\begin{abstract}
This paper presents a new hybrid MC/deterministic method for solving the one-group steady-state
 Boltzmann transport equation based on decomposition of solution in macro and micro components. 
The macro component captures the large scale structure of the solution. 
 It is represented by  angular moments of the high-order transport solution.
The $P_1$ approximation is applied to define the macro component.
The first two angular moments are obtained as a solution of hybrid low-order moment equations with exact closures.
The equation for the micro component is solved using a MC simulation. 
The hybrid two-level system of equations for macro and micro components is solved by fixed-point iteration scheme. 
Numerical results are presented to demonstrate variance reduction of stochastic numerical solution and improvement in computational efficiency.
\end{abstract}
\end{frontmatter}

\section{Introduction}
The particle transport problem is relevant to many applications, so achieving efficient and accurate solutions to this problem is highly desirable. Generally, this problem is solved in two ways: with deterministic methods or with stochastic methods. Deterministic methods solve a discrete form of the high-dimensional transport equation on a finite mesh. This introduces discretization error that converges with
refinement of the mesh.
Stochastic methods, such as Monte Carlo (MC), introduce only stochastic error which decreases with increasing particles. Various hybrid MC/deterministic projection operator schemes using discretized low-order  forms of the transport equation have been shown to improve efficiency and reduce variance \cite{larsen-yang-nse-2008,wolters-nse-2013,Willert-Knoll-Kelley-Park,peterson-morel-ragusa-siam-2016,Lee-Joo-Lee-Smith,vnn-dya-ans-annual-2024,pozulp-haut-brantley-ovilier-JCP-2026,vnn-dya-jctt-2026}. 
There exist also computational methods for multiscale problems applying decomposition of the solution
in  microscopic and macroscopic components for variance reduction \cite{radke-2013}.
This approach has been used to develop hybrid methods for thermal radiative transfer problems \cite{pozulp-haut-brantley-ovilier-JCP-2026}. 

In this paper, we present a new hybrid MC/deterministic method for one-group steady-state 1D transport problems based on the decomposition of the solution in macro and micro components. 
The large-scale structure of the solution is captured by the macro component. This component is formed using angular moments that are the solution of hybrid low-order quasidiffusion (VEF) equations  \cite{goldin-cmmp-1964,auer-mihalas-1970}.
We apply the $P_1$ approximation to define the macro component.
The complementary micro component is solved using a MC simulation. 
The hybrid two-level system of equations for macro and micro components is solved by an iterative process.
We present numerical results demonstrating variance reduction of stochastic numerical solution and improvement in computational efficiency.

\section{Micro-Macro Transport Method}

We consider the steady-state monoenergetic linear transport
problem in 1D slab geometry  defined by
\begin{equation}       \label{eqn:transport}
	\mu\frac{\partial \psi}{\partial x}(x,\mu)+\Sigma_t(x)\psi(x,\mu) \! =  \! 
	\frac{1}{2}\bigg( \Sigma_s(x) \! \!   \int_{-1}^1 \!  \!   \psi(x,\mu') d \mu' +q(x)\bigg) \, ,
\end{equation}
\begin{align} \label{t-eq-bcs}
	\psi(0,\mu) = \psi_{in}^+(\mu)  \ \mbox{for} \ \mu>0,\ \psi(X,\mu) = \psi_{in}^-(\mu) \ \mbox{for} \ \mu<0 \, .
\end{align}
Here $\psi$ is the angular flux,  $x \in [0,X]$ is position, $\mu \in[-1,1]$ is the directional cosine.  $\Sigma_t$ and $\Sigma_s$  are the total and scattering  cross sections, respectively.
$q$ is the external source. $\psi_{in}^{\pm}$ are   prescribed incoming fluxes at boundaries.
The scalar flux of particles is the zeroth angular moment of the angular flux defined by
$\phi(x) = \int_{-1}^{1}  \psi(x,\mu) d \mu$.
The particle current is the first angular moment  of $\psi$ given by
$J(x) = \int_{-1}^{1}  \mu \psi(x,\mu) d \mu$.
We represent the solution of the transport equation (Eq.  \eqref{eqn:transport})  as follows:
\begin{equation} \label{psi-decomp}
	\psi(x,\mu) = \Psi(x,\mu)+\delta\psi(x,\mu) \, ,
\end{equation}
where $\Psi$ is the macro component of the transport solution and $\delta\psi$ is its micro component.
The macro component 
$\Psi$ captures the large-scale behavior of 
the solution represented by the angular moments of $\psi$ and meets the following conditions:
\begin{equation} \label{Psi-conditions}
	\int_{-1}^1\Psi d\mu = \phi,\quad \int_{-1}^1\mu\Psi d\mu=J.
\end{equation}
To formulate the transport problem  for the micro  component $\delta \psi$,
we plug Eq. \eqref{psi-decomp} in Eqs. \eqref{eqn:transport} and  \eqref{t-eq-bcs} to obtain
\begin{equation} \label{eq-d-psi}
	\mu\frac{\partial \delta \psi}{\partial x}(x,\mu)+\Sigma_t(x) \delta\psi(x,\mu)= Q_{\Psi}(x,\mu) \, ,
\end{equation}
\begin{subequations} \label{bcs-d-psi}
	\begin{align}  
		\delta\psi(0,\mu) &= \psi^{+}_{in}(\mu)-\Psi(0,\mu) \  \ \mbox{for} \ \mu>0 \, ,\\
		\delta\psi(X,\mu) &= \psi^{-}_{in}(\mu)-\Psi(X,\mu) \ \ \mbox{for} \ \mu<0 \, ,
	\end{align}
\end{subequations}
where $Q_{\Psi}$ is defined by the macro component as follows: 
\begin{equation}
	Q_{\Psi}(x,\mu) \! =\! \!  \frac{1}{2}\Big(  \Sigma_s(x)\phi(x)+ q(x) \! \Big)    -\mu\frac{\partial \Psi}{\partial x}    (x,\mu)   -\Sigma_t(x)\Psi(x,\mu) \, .
	\label{eqn:HO_source}
\end{equation}
We form the macro component by applying the 
$P_1$-approximation for the angular flux to obtain
\begin{equation} \label{Psi}
	\Psi(x,\mu) = \frac{1}{2} \big(\phi(x) +3 \mu J(x) \big) \, 
\end{equation} 
which 
satisfies the conditions \eqref{Psi-conditions}.
This leads to the right-hand side of the high-order equation for $\delta \psi$  (Eq. \eqref{eq-d-psi}) 
given by
\begin{equation}     \label{Q-Psi}
	Q_{\Psi} =\frac{1}{2} \bigg( \Sigma_s\phi + q    - \mu   \frac{\partial \ }{\partial x} \big( \phi  + 3 \mu J  \big)-  \Sigma_t\big( \phi + 3 \mu J \big) \bigg) \, .
\end{equation}

\clearpage
To define the equations for the scalar flux and current, 
we use
the low-order quasidiffusion (QD) (aka Variable Eddington Factor) equations for these moments given by
\cite{goldin-cmmp-1964,auer-mihalas-1970} 
\begin{subequations} \label{loqd-eqs}
	\begin{equation}         \label{eqn:LOQD_balance}
		\frac{d J}{d x}+(\Sigma_t-\Sigma_s)\phi = q \, ,
	\end{equation}
	\begin{equation}             \label{eqn:LOQD_moment}
		\frac{d}{d x}(E\phi)+\Sigma_t J= 0 \, ,
	\end{equation}
\end{subequations}
\begin{subequations} \label{loqd-bcs}
	\begin{equation}
		J(0) = C_L(\phi(0) - \phi_{in}^+) + J_{in}^+ \, ,   
	\end{equation}  
	\begin{equation}
		J(X) = C_R(\phi(X) - \phi_{in}^-) + J_{in}^- \, ,   
	\end{equation}  
\end{subequations}
where 
\begin{equation}
	\phi_{in}^{\pm} = \pm \int_{0}^{\pm 1} \psi_{in}^{\pm} d \mu \, , \quad
	J_{in}^{\pm} = \pm \int_{0}^{\pm 1} \mu \psi_{in}^{\pm} d \mu \, .
\end{equation}
The exact closures are  defined by 
the QD (Eddington) factor
\begin{equation} \label{qdf}
	E(x) \!  \! = \!\! \int_{-1}^1  \!   \! \mu^2 \Big(\Psi(x,\mu) \! + \! \delta\psi(x,\mu) \Big) d\mu  \bigg/ \!\! \int_{-1}^1  \! \Big(\Psi(x,\mu)  \! + \! \delta\psi(x,\mu) \Big)d\mu
\end{equation}
and the boundary factors given by
\begin{subequations}\label{bcf}
	\begin{equation}
		C_L =  \int_{-1}^0 \mu \Big(\Psi+ \delta\psi \Big) d\mu \bigg/ \int_{-1}^0  \Big(\Psi + \delta\psi \Big) d \mu\bigg|_{x=0} \, ,
	\end{equation}
	\begin{equation} 
		C_R =  \int_0^1 \mu \Big(\Psi  + \delta\psi  \Big) d\mu \bigg/
		\int_0^1  \Big(\Psi  + \delta\psi \Big) d \mu \bigg|_{x=X}\, .
	\end{equation}
\end{subequations}
The closures \eqref{qdf} and \eqref{bcf} are defined with the decomposed angular flux given by Eqs. \eqref{psi-decomp} and \eqref{Psi}.
The micro-macro transport (MMT) method is defined by 
the two-level system of  equations consisting 
of:
(i)  the high-order transport problem   for the micro component $\delta \psi$
(Eqs.	\eqref{eq-d-psi}, \eqref{bcs-d-psi} and \eqref{Q-Psi}) 
and (ii) the low-order QD (LOQD) problem  for $\phi$ and $J$ (Eqs. \eqref{loqd-eqs} and \eqref{loqd-bcs})
which are used to define the macro component $\Psi$ (Eq. \eqref{Psi}).

\section{Spatial Discretization of Equations}

We define a spatial mesh $\{ x_{i}\}_{i=0}^I$.
The cross sections and particle source are piece-wise constant over the set of spatial mesh cells.
The system of the LOQD equations is discretized with a   finite volume scheme
of second-order accuracy
\cite{dya-vyag-ttsp}. 
The balance equation (Eq. \eqref{eqn:LOQD_balance}) 
is integrated over the  $i^{th}$ cell $\omega_i=[x_{i-1}, x_{i}]$.
The first-moment equation (Eq.\eqref{eqn:LOQD_moment}) 
is integrated over the halves of the $i^{th}$ cell. The approximated LOQD equations ($i=1,\ldots,I$) are given by
\begin{subequations}
	\begin{equation}
		J_{i+1}-J_{i}+ \big( \Sigma_{t,i} - \Sigma_{s,i} \big) \Delta x_i \bar\Phi_i = q_i\Delta x_i \, ,  
	\end{equation}
	\begin{equation}
		\bar E_{i} \bar \Phi_{i}  -  E_{i-1}\Phi_{i-1} +   0.5 \Sigma_{t,i} \Delta x_{i} J_{i-1} = 0  \, , 
	\end{equation}
	\begin{equation}
		E_{i}\Phi_{i}  -  \bar E_i \bar \Phi_i  + 0.5 \Sigma_{t,i} \Delta x_{i}J_{i} = 0  \, ,     
	\end{equation}
\end{subequations}
\begin{equation}
	J_{0}=C_L(\Phi_{0}-\phi_{in}^{+})+J_{in}^{+} \, ,
	\quad
	J_{I}=C_R(\Phi_{I}-\phi_{in}^{-})+J_{in}^{-} \, , 
\end{equation}
\begin{equation*}
	\Delta x_i =x_{i}-x_{i-1} \, .
\end{equation*}
The solution of LOQD equations is defined by  
the cell-average scalar fluxes
$\bar \Phi_i~=~\big(\Delta x_i \big)^{-1} \int_{x_{i-1}}^{x_i} \phi(x) d x$,
the   cell-edge scalar fluxes $\Phi_i=\phi(x_i)$, 
and the cell-edge currents $J_i=J(x_i)$.
The coefficients of the LOQD equations are defined by 
the cell-average QD (Eddington) factors $\bar E_i$  and
cell-edge ones $E_i=E(x_i)$. 
The grid functions of the LOQD solution are used to approximate the macro component in space 
and obtain   angular dependent    cell-edge and  cell-average values of $\Psi$ given by
\begin{subequations}\label{Psi-disc}
	\begin{equation} \label{Psi-i}
		\Psi_{i}(\mu) = \frac{1}{2}\bigg(\Phi_{i} +3\mu J_{i }\bigg) \, ,
	\end{equation}
	\begin{equation} \label{bar-Psi-i}
		\bar	\Psi_i(\mu) = \frac{1}{2}\bigg(\bar \Phi_i +\frac{3}{2}\mu (J_{i}+J_{i-1})\bigg) \,  ,
	\end{equation}
\end{subequations}
respectively. Using Eq. \eqref{Psi-disc}, the source term of  the  high-order transport equation  for $\delta \psi$ is approximated
in space by the cell-average values as follows:
\begin{equation}
	\bar Q_{\Psi,i}(\mu) \! = \! \frac{1}{2}\Sigma_{s,i} \bar \Phi_i \! + \! \frac{1}{2}q_i  - \!  \frac{\mu}{\Delta x_i} \Big(   \Psi_{i}(\mu)   -   \Psi_{i-1}(\mu)   \Big) - \! \Sigma_{t,i} \bar \Psi_i (\mu) \, .
	\label{eqn:HO_source_disc}
\end{equation}

\section{Hybrid MMT Algorithm}

The system of equations of the MMT  method is solved by a hybrid MC/deterministic iteration algorithm. 
The transport equation for $\delta \psi$ 
is solved by MC  with  $Q_{\Psi}$ defined by the LOQD solution.
The MC simulations generate quantities needed to calculate 
grid functions of
$E$, $C_L$, $C_R$ 
defined by the micro component $\delta \psi$ and  the macro component $\Psi$.
The angular moments are computed by  means of spatially approximated LOQD equations
with stochastic closures. 
Algorithm~\ref{alg:iteration} presents iterative scheme of the hybrid MMT (HMMT) method, where $\ell$ is the iteration index. 
The  elements of the algorithm are  explained below.
\begin{algorithm}
	\DontPrintSemicolon
	$\ell=0$, $E^{(0)}=\frac{1}{3}$, $C_L^{(0)}=-0.5$, $C_R^{(0)}=0.5$\;
	\While{ $\|  \Phi^{(\ell)} -  \Phi^{(\ell-1)} \|_2  > \varepsilon \|\Phi^{(\ell-1)}\|_2$} {
		$\ell = \ell + 1$\;
		solve hybrid LOQD equations for $\Phi^{(\ell)} , J^{(\ell)} $\;
		solve the transport problem for $\delta \psi  ^{(\ell)}  $ by MC
		to compute  $\delta \varphi_{n}^{(\ell)}$ and
		functionals $E^{(\ell)}, C_L^{(\ell)}, C_R^{(\ell)}$\;
		compute the transport solution 	$\varphi ^{(\ell)}= \Phi^{(\ell)} + \delta   \varphi_{0}^{(\ell)}$\;
	}
	\caption{HMMT Iteration  Scheme
		\label{alg:iteration}} 
\end{algorithm}
Using known closures from the previous iteration, $E$, $C_L$, $C_R$, the LOQD equations are solved for $\Phi$ and $J$.  The MC algorithm then solves the high-order transport problem for the micro component 
with the 
distributed source $\bar Q_{\Psi, i}$ and boundary source terms
\begin{equation}\label{bcs-q}
	Q_{b}^+=-\Psi_{0}(\mu)\ \mbox{for} \ \mu>0, \quad  
	Q_{b}^-=-\Psi_{I}(\mu)\ \mbox{for} \ \mu<0 .
\end{equation} 
defined by Eqs. \eqref{Psi-disc} and \eqref{eqn:HO_source_disc}.
The MC problem does not involve scattering since there is no scattering term in Eq. \eqref{eq-d-psi}.
The form of $\bar{Q}_{\Psi,i}$ and $Q_{b}^{\pm}$ leads to negative sources. The MC algorithm for solving the transport problem for $\delta \psi$ 
uses  particles with positive and negative weights. The sampling of  particles   is performed with stratification. In each spatial cell, a fixed number of particles is sampled with weights corresponding to $\bar{Q}_{\Psi,i}$. For $Q_{b}^{\pm}$, a fixed number of particles is sampled from each boundary. 
We use correlated Monte Carlo sampling \cite{ghoos-et_al-JCP-2016}.
On each iteration, the MC problem for $\delta \psi$ is solved with the same seed. This leads to identical particle histories on each $\ell^{th}$ iteration, just with different weights.
As a result, the statistical noise is fixed during iterations.
The MC simulations compute
\begin{equation}
	\delta \bar \varphi_{n,i} = \frac{1}{\Delta x_i} \int_{\omega_i} \int_{-1}^{1} \mu^{n} \delta \psi d \mu d x \, , \quad n=0,2 \,
\end{equation}
that is, the zeroth and second moments of the micro component $\delta\psi$  in  mesh cells $\omega_i$
and 
\begin{equation}
	\delta  \varphi_{n,i} =  \int_{-1}^{1} \mu^{n} \delta \psi d \mu  \Big|_{x=x_i} \, , \quad n=0,2 \, .
\end{equation}
The QD (Eddington) factors are calculated using the MC solution of the micro-component problem and   the solution of  LOQD problem as follows:
\begin{equation}
	\bar E_i = \frac{\frac{2}{3} \bar \Phi_i + \delta \bar \varphi_{2,i}}{ 2\bar \Phi_i  + \delta \bar \varphi_{0,i} } \, , \quad
	E_i = \frac{\frac{2}{3}   \Phi_i + \delta   \varphi_{2,i}}{  2\Phi_i +  \delta  \varphi_{0,i} } \, .
\end{equation}
The boundary factors $C_L$ and $C_R$ (Eq. \eqref{bcf}) are computed similarly using corresponding MC tallies. 
Finally, the solution of the transport problem  \eqref{eqn:transport} and \eqref{t-eq-bcs}  is the
particle scalar flux defined according to  the micro-macro decomposition (Eq. \eqref{psi-decomp})
and given by
\begin{equation}
	\bar	\varphi_i = \bar \Phi_i + \delta \bar \varphi_{0,i}  \, , \quad  \ 	\varphi_i = \Phi_i + \delta   \varphi_{0,i}   \,  .
\end{equation}
Here  $\bar	\varphi_i$ and $\varphi_i $ are the cell-average and cell-edge values, respectively.

\section{Numerical Results}

\begin{figure}[h]
\centering
\includegraphics[scale=1.1]{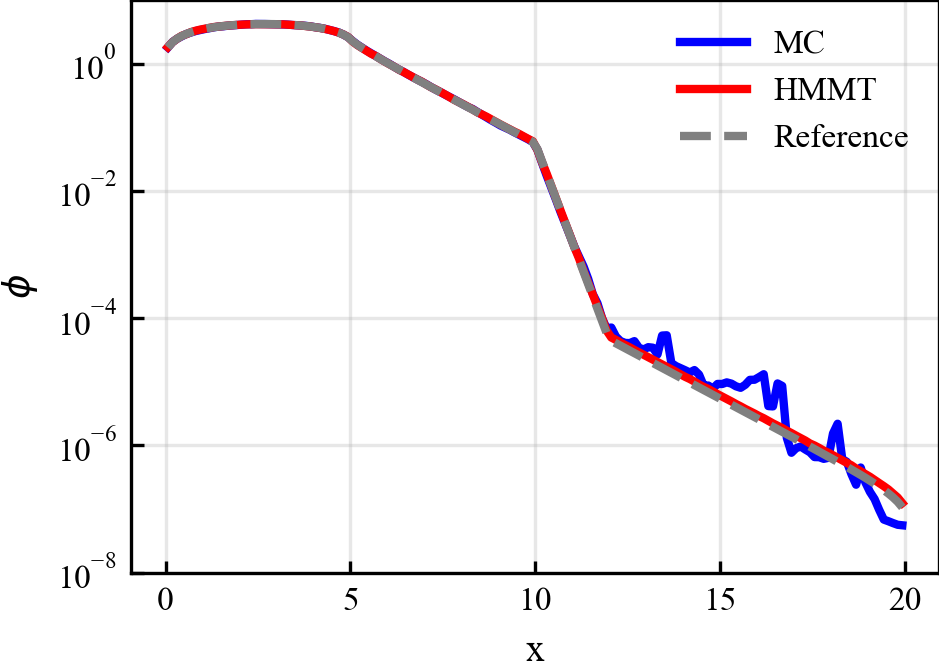}
\caption{Batch mean  of  the scalar flux}
\label{fig:flux}
\end{figure}

\begin{figure}[h]
\centering
\includegraphics[scale=1.1]{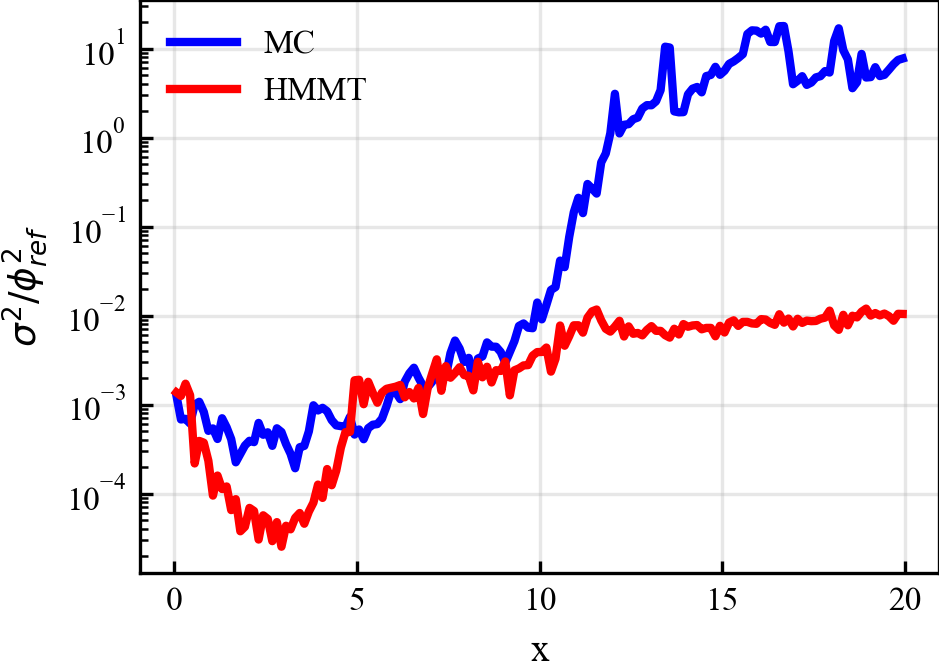}
\caption{Batch mean relative variance in the scalar flux}
\label{fig:rel_var}
\end{figure}

We present numerical results for a  test problem consisting of a heterogeneous slab with two materials and three regions driven by a distributed source. The details are in Table~\ref{tab:problem_description}. The region in $[10,12]$ is an optically thick barrier. A reference solution to this problem is available in Ref. \cite{Novellino2026-diss}.
\begin{table}[h]
	\caption{Test Problem}
	\centering
	\begin{tabular}{|c|c|c|c|c|}
		\hline 
		$x $& $[0,5]$ & $[5,10]$ & $[10,12]$& $[12,20]$  \\ 	\hline 
		$\Sigma_t$ & $1$ & $1$ & $5$& $1$ \\ 	\hline 
		$\Sigma_s$ & $0.8$ & $0.8$ & $4$& $0.8$ \\	\hline 
		$q$ & $1$ & $0$ & $0$& $0$ 	 \\ \hline
	\end{tabular}
	\label{tab:problem_description}
\end{table}
To  assess performance of the HMMT method,
we compute the mean and variance of the numerical solution  on  a uniform spatial mesh with   $\Delta x =2^{-3}$
using 20 batches with different random number seeds but the same number of particles, $N=2^{13}$. 
For all results in this work, a fixed convergence criterion of $\varepsilon=10^{-3}$ was used.
For all the cases we ran, the highest number of iterations needed to converge was 15.
Figure \ref{fig:flux} presents the batch mean of the scalar flux   calculated by traditional MC and HMMT. The MC results were computed using implicit capture (IC) for variance reduction. These results show that the solution of the HMMT method is more accurate in the subdomain $x\in[12,20]$. The analysis also showed that the HMMT solution is more accurate in the source region. 
The $L_2$-norm of the relative error in the HMMT and MC solutions is 0.48 and 1.3, respectively.
In Figure \ref{fig:rel_var} is the batch relative variance of the numerical solutions.
These results demonstrate that the HMMT method reduced variance in the source region as well as in the subdomain containing the barrier and behind it.
\FloatBarrier

\begin{figure}[h]
\centering
\includegraphics[scale=1.1]{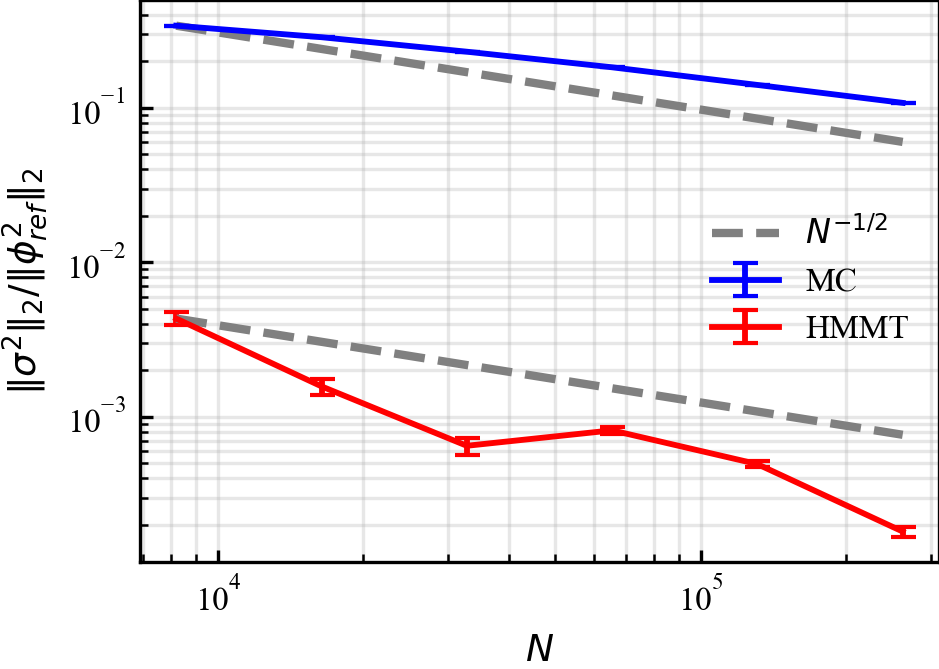}
\caption{Relative variance convergence with increasing $N$}
\label{fig:var_convergence}
\end{figure}

\begin{figure}[h]
\centering
\includegraphics[scale=1.1]{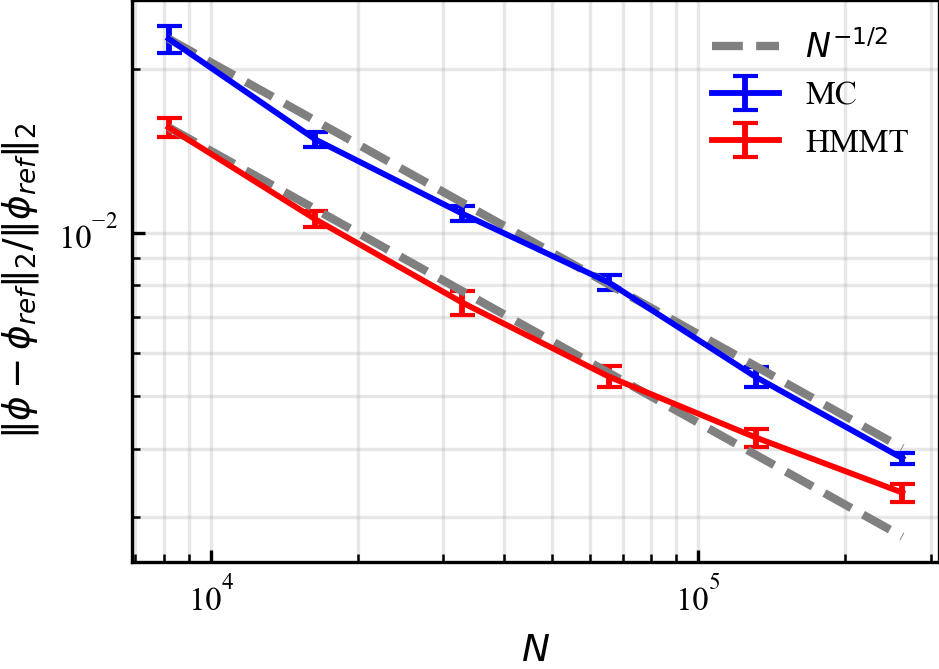}
\caption{Relative error convergence with increasing $N$}
\label{fig:error_convergence}
\end{figure}

\begin{figure}[h!]
\centering
\includegraphics[scale=1.1]{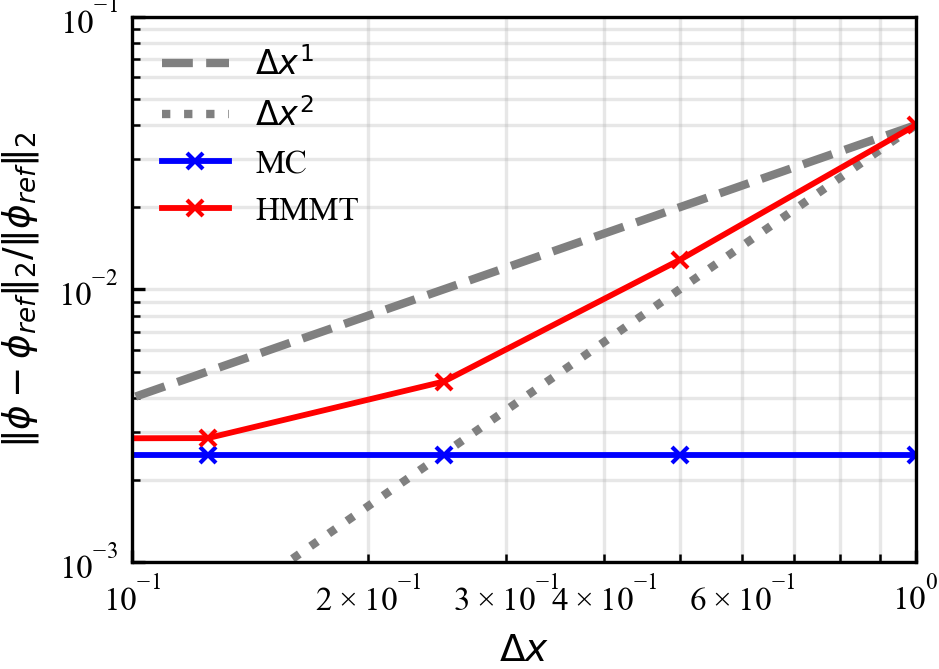}
\caption{Relative error convergence with decreasing $\Delta x$}\label{fig:error_convergence_dx}
\end{figure}

\begin{figure}[h!]
	\centering
	\includegraphics[]{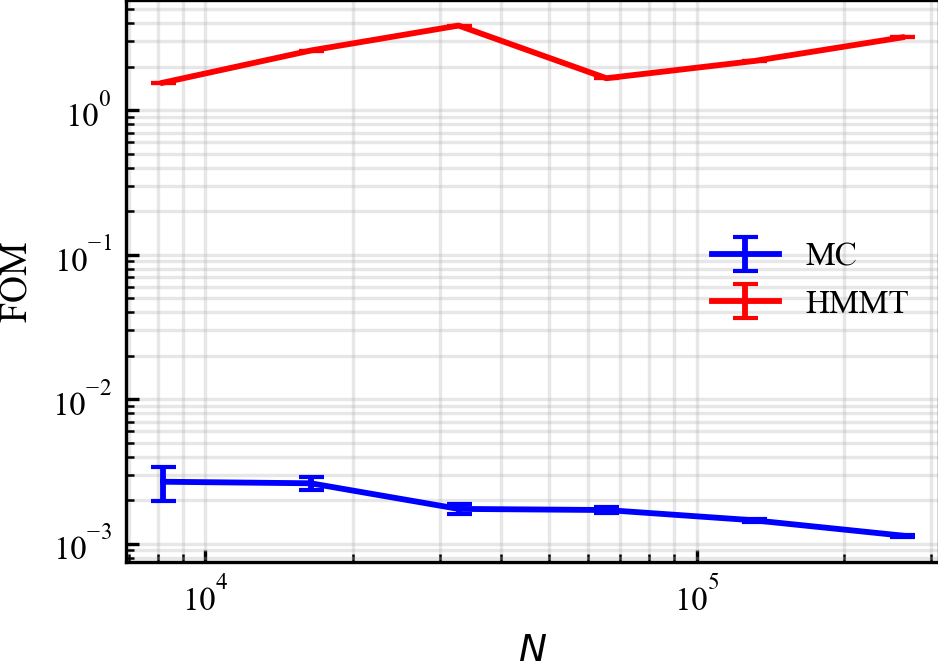}
	\caption{FOM convergence with increasing $N$}
	\label{fig:FOM_convergence}
\end{figure}

The error in the HMMT solution has two parts: discretization error which decreases
as $\Delta x \to 0$,
and stochastic error which decreases 
as $N \to \infty$.
The HMMT method  has theoretical convergence  
 $\mathcal{O}(\Delta x)~+~\mathcal{O}(N^{-1/2})$. We present results for convergence studies in both $\Delta x$ and $N$. Beginning with a fixed $\Delta x = 2^{-3}$ we vary the number of histories over the range $N=[2^{13},2^{14},2^{15},2^{16},2^{17},2^{18}]$. 	 We run 20 batches for each $N$ and compute mean errors and variances.
The convergence of the relative variance and error in the 2-norm are plotted in Figures \ref{fig:var_convergence} and \ref{fig:error_convergence}, respectively. Figure \ref{fig:var_convergence} shows that HMMT solutions have about 2 orders of magnitude less variance in the 2-norm compared to MC with IC. Figure \ref{fig:error_convergence} demonstrates that the methods achieve near the theoretical $N^{-1/2}$ convergence, with the  HMMT solutions having greater accuracy.
We note that the convergence of the HMMT solution deviates from the $N^{-1/2}$  at the highest number of histories, $N=2^{18}$. 
This is the effect of discretization error which becomes comparable in magnitude with the stochastic error for this combination of $N$ and $\Delta x$.
We now fix the number of histories at $N=2^{18}$ and vary $\Delta x =[1,2^{-1},2^{-2},2^{-3},2^{-4},2^{-5},2^{-6}]$. The 2-norm relative errors for these cases are shown in Figure  \ref{fig:error_convergence_dx}. 	
The relative error of the MC solution is constant since it has no discretization error. 
For coarse meshes, the HMMT error converges with greater than first-order rate.
After $\Delta x=2^{-3}$ the stochastic error is dominant.
To measure efficiency of the algorithms, we consider the figure of merit (FOM)  defined as	$FOM=(\|\sigma^2\|_2\tau)^{-1}$, where $\tau$ is the runtime. 
This FOM should remain constant over various $N$. In Figure \ref{fig:FOM_convergence} the FOM is plotted for the range of histories
for the spatial mesh with $\Delta x =2^{-3}$.
The FOM for both MC and HMMT remains 
approximately
constant with increasing number of histories.
These results indicate that
the FOM achieved with the HMMT method is over 100 times greater than with MC. 
The FOM improvement can be partially attributed to reduced runtimes with the HMMT method. 
The mean runtimes for the scaling study are shown in Table \ref{tab:runtimes}. 
Reduced runtime is  possible  since the MC portion of the HMMT method has no scattering events which often dominate runtime.
\begin{table}[h]
\centering
\caption{Runtimes  [s] versus $N$ for $\Delta x = 2^{-3}$   \label{tab:runtimes}}
\begin{tabular}{|c|c|c|}
	\hline
$N$ & HMMT& MC   \\
\hline
$2^{13}$ & 1.40e+02 & 1.46e+02 \\ \hline
$2^{14}$ & 2.37e+02 & 2.92e+02 \\ \hline
$2^{15}$ & 4.02e+02 & 5.84e+02 \\ \hline
$2^{16}$ & 7.69e+02 & 1.17e+03 \\ \hline
$2^{17}$ & 9.18e+02 & 2.34e+03 \\ \hline
\end{tabular}
\end{table}

\section{Conclusion}

We have presented	 a new hybrid method
for solving   transport problems in 1D slab geometry. It is based on a decomposition of the angular flux in macro and micro components. 
The numerical results show that the proposed HMMT method reduces variance and significantly improves figure of merit compared to the MC algorithm with IC.
Future research will include  detailed analysis of the HMMT method on various test problems, development of different versions of the method
using various definitions of the macro component, application of other types of low-order equations, and advanced discretization schemes. 

\section{Acknowledgments}

This work was supported by the Center for Advancing the Radiation
Resilience of Electronics (CARRE),
a PSAAP-IV project funded by
the Department of Energy  
National Nuclear Security Administration, award number: DE-NA0004268.

\bibliographystyle{elsarticle-num}
\bibliography{cas-dya-tans-2026-arxiv}

\begin{thebibliography}{10}
\expandafter\ifx\csname url\endcsname\relax
  \def\url#1{\texttt{#1}}\fi
\expandafter\ifx\csname urlprefix\endcsname\relax\def\urlprefix{URL }\fi
\expandafter\ifx\csname href\endcsname\relax
  \def\href#1#2{#2} \def\path#1{#1}\fi

\bibitem{larsen-yang-nse-2008}
{E. W. Larsen}, {J. Yang}, A functional {M}onte {C}arlo method for k-eigenvalue
  problems, Nuclear Science and Engineering 159 (2008) 107--126.

\bibitem{wolters-nse-2013}
{E. R. Wolters}, {E. W. Larsen}, {W. R. Martin}, Hybrid {M}onte
  {C}arlo–{CMFD} methods for accelerating fission source convergence, Nuclear
  Science and Engineering 174 (2013) 286–299.

\bibitem{Willert-Knoll-Kelley-Park}
J.~Willert, C.~T. Kelley, D.~A. Knoll, Hybrid deterministic/{M}onte {C}arlo
  neutronics, SIAM Journal on Scientific Computing 35~(5) (2013) S62--S83.

\bibitem{peterson-morel-ragusa-siam-2016}
J.~R. Peterson, J.~E. Morel, J.~C. Ragusa, Residual {M}onte {C}arlo for the
  one-dimensional particle transport equation, SIAM Journal on Scientific
  Computing 38~(6) (2016) B941--B961.

\bibitem{Lee-Joo-Lee-Smith}
M.~J. Lee, H.~G. Joo, D.~Lee, K.~Smith, Coarse mesh finite difference
  formulation for accelerated {M}onte {C}arlo eigenvalue calculation, Annals of
  Nuclear Energy 65 (2014) 101--113.

\bibitem{vnn-dya-ans-annual-2024}
V.~N. Novellino, D.~Y. Anistratov, Analysis of hybrid {MC}/deterministic
  methods for transport problems based on low-order equations discretized by
  finite volume schemes, Transactions of American Nuclear Society 130 (2024)
  408--411.

\bibitem{pozulp-haut-brantley-ovilier-JCP-2026}
M.~Pozulp, T.~Haut, P.~Brantley, S.~Olivier, A hybrid {M}onte
  {C}arlo-deterministic second moment method with efficient variance reduction,
  Journal of Computational Physics 565 (2026) 115121.

\bibitem{vnn-dya-jctt-2026}
V.~N. Novellino, D.~Y. Anistratov, Multi-level hybrid {M}onte
  {C}arlo/deterministic methods for particle transport problems, Journal of
  Computational and Theoretical Transport (2026) 1--34.

\bibitem{radke-2013}
G.~A. Radtke, J.-P.~M. Péraud, N.~G. Hadjiconstantinou, On efficient
  simulations of multiscale kinetic transport, Philosophical Transactions of
  the Royal Society A: Mathematical, Physical and Engineering Sciences
  371~(1982) (2013) 20120182.

\bibitem{goldin-cmmp-1964}
V.~Y. Gol'din, A quasi-diffusion method of solving the kinetic equation, Comp.
  Math. and Math. Phys. 4 (1964) 136--149.

\bibitem{auer-mihalas-1970}
L.~H. Auer, D.~Mihalas, On the use of variable {E}ddington factors in non-{LTE}
  stellar atmospheres computations, Monthly Notices of the Royal Astronomical
  Society 149 (1970) 65--74.

\bibitem{dya-vyag-ttsp}
D.~Y. Anistratov, V.~Y. Gol'din, Nonlinear methods for solving particle
  transport problems, Transport Theory and Statistical Physics 22 (1993)
  42--77.

\bibitem{ghoos-et_al-JCP-2016}
K.~Ghoos, W.~Dekeyser, G.~Samaey, P.~Börner, M.~Baelmans, Accuracy and
  convergence of coupled finite-volume / {M}onte {C}arlo codes for plasma edge
  simulations of nuclear fusion reactors, Journal of Computational Physics 322
  (2016) 162--182.

\bibitem{Novellino2026-diss}
V.~N. Novellino, Multilevel hybrid {M}onte {C}arlo / deterministic methods for
  neutral particle transport, Ph.D. thesis, North Carolina State University
  (May 2026).

\end{thebibliography}

\end{document}